# LIMIT LAWS OF ESTIMATORS FOR CRITICAL MULTI-TYPE GALTON–WATSON PROCESSES

By Zhiyi Chi

*University of Chicago*

We consider the asymptotics of various estimators based on a large sample of branching trees from a critical multi-type Galton–Watson process, as the sample size increases to infinity. The asymptotics of additive functions of trees, such as sizes of trees and frequencies of types within trees, a higher-order asymptotic of the "relative frequency" estimator of the left eigenvector of the mean matrix, a higher-order joint asymptotic of the maximum likelihood estimators of the offspring probabilities and the consistency of an estimator of the right eigenvector of the mean matrix, are established.

**1. Introduction.** This article considers the asymptotics of estimators associated with critical multi-type Galton–Watson (GW) processes. A GW process is called critical if the largest eigenvalue of its mean matrix is 1 (see below for details). For such a process, a branching tree is finite with probability 1, but the expectation of its size is infinite. The estimators considered here are based on a large sample of terminating branching trees, and the asymptotics refer to the probabilistic behavior as the sample size $n \to \infty$.

The study on large sample asymptotics of parameter estimators for simple (i.e., single type) GW processes has a quite long history (cf. [20]). The idea of using increasingly large sample of individual trees for estimation dates from as early as [24], and much progress has been made since then (cf. [7, 8] and references therein). This setting of parameter estimation is widely used in computational linguistics [5, 16], where large samples of tree-structured parses of sentences are available. On the general issues of parameter estimation or asymptotics related to simple or multi-type GW processes, there is now extensive literature available (e.g., [1, 3, 6, 10, 12, 14, 15, 17, 19, 21, 22, 23] and references therein). For a multi-type GW process, in addition to the estimation of offspring probabilities associated with different types, there is a unique problem, namely the estimation of the left and right eigenvectors of the mean matrix of the process. Both estimation problems will be dealt with later in the article. Indeed, by the Perron–Frobenius theorem, there is a unique pair of left and right eigenvectors that satisfy certain conditions. We will refer to these two eigenvectors as Frobenius eigenvectors. The estimators considered in the article give consistent estimation for both of them.

Since the estimation relies on the asymptotics of the total size of sample branching trees, we shall first establish results in this regard. For simple GW processes, it is well known that the distribution of the size of a tree, that is, its total progeny, belongs to the domain of attraction of a stable law of exponent $\frac{1}{2}$ [11]. The size of a tree is an additive function of the tree (cf. [13] and Section 4). Under a critical GW process, additive functions exhibit very different asymptotic behavior from those under a subcritical process. The consistency of the estimator for the left Frobenius eigenvector of the mean matrix as well as that for the offspring probabilities is a simple consequence of a general result on the asymptotics of an additive function (cf. Theorem 5). For the left Frobenius eigenvector, the estimator consists of relative frequencies of types. From the offspring probabilities, the estimators are the well-known maximum likelihood estimators, which also take the form of relative frequencies. In analogy to the central limit theorem, the fluctuations of these estimators around their limits are also of interest and characterized with non-Gaussian behavior. The estimation of the Frobenius right eigenvector, on the other hand, follows a completely different approach.

The other sections of the article are organized as follows. The main results are stated in Section 2. Some well-known or standard results are collected in Section 3 for later use. Section 4 demonstrates a general result on the asymptotics of additive functions of sample branching trees. In Sections 5 and 6, some finer asymptotics for the estimators of the left Frobenius eigenvector and the offspring probabilities are studied. Finally, a consistent estimator of the right Frobenius eigenvector is given in Section 7.

In the rest of this section we shall fix the notation for the article. Throughout, we use $V \in \mathbb{N}$ as a generic notation for the number of particle types in a GW process. Without loss of generality, let the set of types be $\mathcal{V} = \{1, \ldots, V\}$. For simplicity, the topology of a sample branching tree will be ignored, and a branching rule in which a particle of type $k$ generates $n_1$ offspring of type 1, $n_2$ offspring of type 2 and so on is denoted by $k \to \mathbf{n}$, where $\mathbf{n} = (n_1, \ldots, n_V) \in \mathcal{V}^* := (\{0\} \cup \mathbb{N})^V$. Indeed, when the topology of a tree needs to be taken into account, one can denote by $\mathbf{n}$ a finite string consisting of elements in $\mathcal{V}$, and by $\mathcal{V}^*$ the set of all such strings, $n_s$ is still the number of particles of type $s$ in $\mathbf{n}$ and results established in this article still hold.

Given $k \in \mathcal{V}$, the offspring probability distribution on $\mathcal{V}^*$, namely the probability distribution on $k \to \mathbf{n}$, is denoted by $p_k(\mathbf{n})$. Let $P_k$ denote the probability distribution determined by the offspring probability distributions $p_s$, $s \in \mathcal{V}$, on branching trees rooted with a particle of type $k$, and let $E_k$ denote the expectation under $P_k$. For each sample tree $\omega$, $k \in \mathcal{V}$ and $\mathbf{n} \in \mathcal{V}^*$, denote

$$|\omega| = \text{number of particles in } \omega,$$

$$f(k;\omega) = \text{number of particles of type } k \text{ in } \omega,$$

$$\mathbf{f}(\omega) = (f(1;\omega), \ldots, f(V;\omega)),$$

$$f(k \to \mathbf{n}; \omega) = \text{number of times the rule } (k \to \mathbf{n}) \text{ is applied in } \omega,$$

$$|\mathbf{n}| = \sum_{s \in \mathcal{V}} n_s.$$

Then for any finite sample tree $\omega$ rooted with a particle of type $k$,

$$P_k(\omega) = \prod_{s \in \mathcal{V}} \prod_{\mathbf{n} \in \mathcal{V}^*} p_s(\mathbf{n})^{f(s \to \mathbf{n}; \omega)}, \qquad k \in \mathcal{V}.$$

The mean matrix of a GW process is a $V \times V$ matrix $M$, with its $(k,l)$th entry

$$M(k,l) = M_{kl} = \sum_{\mathbf{n} \in \mathcal{V}^*} p_k(\mathbf{n}) n_l \geq 0.$$

In general, a matrix $M$ is called nonnegative (resp. positive) if all its entries are nonnegative (resp. positive). When $M$ is a square matrix, it is called irreducible if $M^p$ is positive for some $p \in \mathbb{N}$.

Henceforth, any $\mathbf{v} \in \mathbb{C}^V$ will be regarded as a row vector, and its transpose $\mathbf{v}^t$ as a column vector. For two vectors $\mathbf{v}$ and $\mathbf{u}$, denote $\mathbf{v} \cdot \mathbf{u} = \sum_{s \in \mathcal{V}} v_s u_s$. For any scalar $a$, denote $\mathbf{a} = (a, \ldots, a)$.

Finally, recall the following fundamental result (cf. [2], page 185, and [16]).

PERRON–FROBENIUS THEOREM. *Let $M$ be a nonnegative matrix indexed by $\mathcal{V} \times \mathcal{V}$. Then $M$ has an eigenvalue $\lambda \in [0, \infty)$ such that no other eigenvalue of $M$ has absolute value greater than $\lambda$; and there are nonnegative vectors $\mathbf{v}, \mathbf{u} \in \mathbb{R}^V$ satisfying $\lambda \mathbf{v} = \mathbf{v} M$, $\lambda \mathbf{u}^t = M \mathbf{u}^t$. Moreover, if $M$ is irreducible, then $\lambda$ is a simple eigenvalue, and $\mathbf{v}$ and $\mathbf{u}$ are positive and can be chosen in such a way that*

$$(1.1) \qquad \sum_{s \in \mathcal{V}} v_s = 1, \qquad \mathbf{v} \cdot \mathbf{u} = 1, \qquad M^n = \lambda^n \mathbf{u}^t \mathbf{v} + R_n,$$

*where $\max_{k,l \in \mathcal{V}} |R_n(k,l)| = O(\alpha^n)$ with $0 \leq \alpha < \lambda$. Indeed, for all $n \geq 1$, $R_n = R_1^n$ and $\sigma(R_1) := \sup_{|\mathbf{x}|=1} |R_1 \mathbf{x}^t| < \lambda$.*

It is easy to see that eigenvectors satisfying (1.1) are unique. We will refer to them as Frobenius eigenvectors and denote them by $\mathbf{v}$ and $\mathbf{u}$, respectively.

**2. Main results.** Define measure $Q$ on $\mathcal{V}^*$ by
$$Q(\mathbf{n}) = \sum_{s \in \mathcal{V}} v_s p_s(\mathbf{n}).$$

It is not hard to see $Q$ is a probability measure. Henceforth, we will denote by $\mathbf{X}$ the identity function on $\mathcal{V}^*$ [i.e., $\mathbf{X}(\mathbf{n}) = \mathbf{n}$] and assume

(2.1) $$E_Q(\mathbf{X} \cdot \mathbf{u})^2 < \infty.$$

Given $k \in \mathcal{V}$, let $\omega_1, \omega_2, \ldots$ be i.i.d. trees sampled from $P_k$. For the asymptotics of the total numbers of particles of different types in the trees, convergence to a joint stable distribution can be established, which generalizes the well-known result that the distribution of the total progeny of a critical simple GW process is in the domain of attraction of a stable law with exponent $\frac{1}{2}$ (cf. [11], Theorem 13.1, and [4], Theorem 9.34).

THEOREM 1. *Let $\mathbf{v}$ and $\mathbf{u}$ be the Frobenius eigenvectors given in (1.1), with $\lambda = 1$. Suppose (2.1) holds. Then*

(2.2) $$\frac{1}{N^2} \sum_{n=1}^{N} \mathbf{f}(\omega_n) \xrightarrow{\mathcal{D}} \frac{\mathbf{v} u_k^2 \xi}{H(\mathbf{u})}, \qquad N \to \infty,$$

*where $\xi$ is a real-valued stable random vector with characteristic function*

(2.3) $$E[e^{it\xi}] = \exp\{-(1 - i\,\mathrm{sign}(t))\sqrt{|t|}\,\}$$

*and $H$ is defined on $\mathbb{C}^V$ by*

(2.4) $$H(\mathbf{z}) = E_Q(\mathbf{X} \cdot \mathbf{z})^2 - \sum_{s \in \mathcal{V}} v_s z_s^2 \qquad \forall\, \mathbf{z} = (z_1, \ldots, z_V) \in \mathbb{C}^V.$$

*Note that $H(\mathbf{u}) > 0$.*

Theorem 1 is a special case of the asymptotics of functions of the form ([13], page 167).
$$G(\omega) = \sum_{s \in \mathcal{V}} \sum_{\mathbf{n} \in \mathcal{V}^*} g_s(\mathbf{n}) f(s \to \mathbf{n}; \omega),$$

which we will refer to as "additive functions." There are many choices for $g_s$. For example, if $g_s(\mathbf{n}) = \mathbf{1}_{\{|\mathbf{n}|=0\}}$, then $G(\omega)$ is equal to the total number of terminals in $\omega$. In Section 4, it will be shown that, under suitable conditions, $\frac{1}{N^2} \sum_{n=1}^{N} G(\omega_n)$ converges in distribution, and Theorem 1 immediately follows.

From (2.2) and the equivalence between convergence in distribution to a constant and convergence in probability to the same constant, it follows that

(2.5) $$\hat{\mathbf{v}}_N = \hat{\mathbf{v}}_N(\omega_1, \ldots, \omega_N) = \frac{\sum_{n=1}^{N} \mathbf{f}(\omega_n)}{\sum_{n=1}^{N} |\omega_n|} \xrightarrow{P} \mathbf{v}, \qquad \omega_1, \omega_2, \ldots \text{ i.i.d.} \sim P_k.$$

To find finer asymptotics of $\hat{\mathbf{v}}_N$, we next consider the limit of the characteristic functions of $N^\alpha(\hat{\mathbf{v}}_N - \mathbf{v})$, as $N \to \infty$. It turns out that $\alpha = 1$ is the correct scaling. Since

$$N(\hat{\mathbf{v}}_N - \mathbf{v}) = \frac{(1/N)\sum_{n=1}^N (\mathbf{f}(\omega_n) - \mathbf{v}|\omega_n|)}{(1/N^2)\sum_{n=1}^N |\omega_n|}$$

and

$$\frac{1}{N^2}\sum_{n=1}^N |\omega_n| = \frac{1}{N^2}\sum_{n=1}^N \mathbf{1} \cdot \mathbf{f}(\omega_n) \xrightarrow{\mathcal{D}} \frac{u_k^2 \xi}{H(\mathbf{u})},$$

instead of directly dealing with $N(\hat{\mathbf{v}}_N - \mathbf{v})$, we shall consider the limit of the joint characteristic functions of random vectors

$$\left(\frac{1}{N}\sum_{n=1}^N (\mathbf{f}(\omega_n) - \mathbf{v}|\omega_n|), \frac{1}{N^2}\sum_{n=1}^N |\omega_n|\right).$$

Recall the matrices $R_n$ in (1.1). Since the GW processes are critical, $R_n \to 0$ at an exponential rate. Define matrix

(2.6) $$\Lambda = \left(\sum_{n=0}^\infty R_n\right)(I - \mathbf{1}^t\mathbf{v}) = (I - R)^{-1}(I - \mathbf{1}^t\mathbf{v}),$$

where $I$ is the $V \times V$ identity matrix. Then $\Lambda \neq 0$, as is seen from $\mathbf{1}^t\mathbf{v} \neq I$ and

(2.7)
$$(M - I)\Lambda = \sum_{n=0}^\infty (M - I)(M^n - \mathbf{u}^t\mathbf{v})(I - \mathbf{1}^t\mathbf{v})$$
$$= \sum_{n=0}^\infty (M^{n+1} - M^n)(I - \mathbf{1}^t\mathbf{v}) = \mathbf{1}^t\mathbf{v} - I.$$

THEOREM 2. *Under the same conditions as in Theorem 1,*

$$\left(\frac{1}{N}\sum_{n=1}^N (\mathbf{f}(\omega_n) - \mathbf{v}|\omega_n|), \frac{1}{N^2}\sum_{n=1}^N |\omega_n|\right) \xrightarrow{\mathcal{D}} (\mathbf{Z}, W), \qquad N \to \infty,$$

*where* $(\mathbf{Z}, W) \in \mathbb{R}^V \times \mathbb{R}$ *has characteristic function*

$$E[\exp(i\mathbf{c} \cdot \mathbf{Z} + iKW)] = e^{z(\mathbf{c},K)u_k + i\eta_k}, \qquad \mathbf{c} \in \mathbb{R}^V, \ K \in \mathbb{R},$$

*such that* $\eta_k$ *is the kth component of* $\boldsymbol{\eta} \in \mathbb{R}^V$ *given by*

(2.8) $$\boldsymbol{\eta}^t = \Lambda \mathbf{c}^t,$$

*and* $z(\mathbf{c}, K)$ *the (unique) solution with negative real part to*

(2.9) $$z^2 + \frac{2Ai}{H(\mathbf{u})}z + \frac{1}{H(\mathbf{u})}(B + 2Ki) = 0,$$

*where*

$$A = \operatorname{Cov}_Q(\mathbf{X} \cdot \mathbf{u}, \mathbf{X} \cdot \boldsymbol{\eta}) - \sum_{s \in \mathcal{V}} v_s u_s(\eta_s - c_s) - \mathbf{c} \cdot \mathbf{v}, \tag{2.10}$$

$$B = -\operatorname{Var}_Q(\mathbf{X} \cdot \boldsymbol{\eta}) + \sum_{s \in \mathcal{V}} v_s(c_s - \eta_s)^2 - (\mathbf{c} \cdot \mathbf{v})^2. \tag{2.11}$$

*Then immediately one gets* $N(\hat{\mathbf{v}}_N - \mathbf{v}) \xrightarrow{\mathcal{D}} \frac{\mathbf{Z}}{W}$, *as* $N \to \infty$.

The vector $\hat{\mathbf{v}}_N$ consists of the relative frequencies of types in the population of particles in $\omega_1, \ldots, \omega_N$. Likewise, we can consider the relative frequencies of branching rules in $\omega_1, \ldots, \omega_N$. Fix $j \in \mathcal{V}$ and $\mathbf{n} \in \mathcal{V}^*$. For $\omega_1$, $\omega_2, \ldots$ i.i.d. $\sim P_k$, define

$$\hat{p}_{j,N}(\mathbf{n}) = \hat{p}_{j,N}(\mathbf{n}; \omega_1, \ldots, \omega_N) = \frac{\sum_{n=1}^{N} f(j \to \mathbf{n}; \omega_n)}{\sum_{n=1}^{N} f(j; \omega_n)}.$$

From Corollary 1 in Section 4, it is seen that $\hat{p}_{j,N}$ is consistent. That is, for $\omega_1, \omega_2, \ldots$ i.i.d. $\sim P_k$, $\hat{p}_{j,N}(\mathbf{n}) \xrightarrow{P} p_j(\mathbf{n})$, as $N \to \infty$. To get finer asymptotics of $\hat{p}_{j,N}(\mathbf{n})$, following Theorem 2, consider the limit of the joint characteristic functions of

$$\left( \frac{1}{N} \sum_{n=1}^{N} (f(j \to \mathbf{n}; \omega_n) - p_j(\mathbf{n}) f(j; \omega_n)), \mathbf{n} \in \mathcal{V}^*; \frac{1}{N^2} \sum_{n=1}^{N} f(j; \omega_n) \right).$$

Because $\mathcal{V}^*$ may have infinitely many elements, to avoid potential difficulty, we only consider the joint asymptotic of a finite number of relative frequencies of branching rules $(j \to \mathbf{n})$.

THEOREM 3.  *Assume the same conditions as in Theorem* 1. *Given* $j \in \mathcal{V}$ *and* $\mathbf{n}_1, \ldots, \mathbf{n}_M \in \mathcal{V}^*$, *let* $\mathbf{F}(\omega) = (f(j \to \mathbf{n}_1), \ldots, f(j \to \mathbf{n}_M))$, $\mathbf{q} = (p_j(\mathbf{n}_1), \ldots, p_j(\mathbf{n}_M))$. *Then*

$$\left( \frac{1}{N} \sum_{n=1}^{N} (\mathbf{F}(\omega_n) - \mathbf{q} f(j; \omega_n)), \frac{1}{N^2} \sum_{n=1}^{N} f(j; \omega_n) \right) \xrightarrow{\mathcal{D}} (\mathbf{Z}, W), \qquad N \to \infty,$$

*where* $(\mathbf{Z}, W) \in \mathbb{R}^M \times \mathbb{R}$ *has characteristic function*

$$E[\exp(i\mathbf{c} \cdot \mathbf{Z} + iKW)] = e^{z(\mathbf{c}, K) u_k}, \qquad \mathbf{c} = (c_1, \ldots, c_M) \in \mathbb{R}^M, \ K \in \mathbb{R},$$

*such that* $z(\mathbf{c}, K)$ *is the (unique) solution with negative real part to*

$$z^2 + \frac{2 v_j A i}{H(\mathbf{u})} z + \frac{v_j}{H(\mathbf{u})} (B + 2Ki) = 0, \tag{2.12}$$

*where*

$$(2.13) \quad A = \sum_{\mu=1}^{M} c_\mu p_j(\mathbf{n}_\mu)(\mathbf{n}_\mu \cdot \mathbf{u}) - (\mathbf{c} \cdot \mathbf{q})u_j,$$

$$(2.14) \quad B = \sum_{\mu=1}^{M} p_j(\mathbf{n}_\mu)c_\mu^2 - (\mathbf{c} \cdot \mathbf{q})^2.$$

*Then immediately one gets $N(\hat{p}_{j,N}(\mathbf{n}) - p_j(\mathbf{n})) \xrightarrow{\mathcal{D}} \frac{\mathbf{Z}}{W}$, as $N \to \infty$.*

Formula (2.5) gives an estimator of the left Frobenius eigenvector $\mathbf{v}$ of the mean matrix $M$. While the right Frobenius eigenvector $\mathbf{u}$ of $M$ occurs in the asymptotics of the estimator for $\mathbf{v}$, it is clear how to use relative frequencies to directly estimate $\mathbf{u}$. So we consider an alternative approach to the estimation of $\mathbf{u}$. Given a tree $\omega$, for each node $x \in \omega$, let $|x|$ denote its "depth," that is, the number of edges on the shortest path from $x$ to the root of $\omega$. Denote

$$S(\omega, \lambda) = \sum_{x \in \omega} \lambda^{|x|}$$

whenever the sum on the right-hand side is well defined. Recall that $\mathbf{v}$ and $\mathbf{u}$ denote the positive left and right eigenvectors of $M$, respectively, such that the sum of the components of $\mathbf{v}$ is equal to 1, and $\mathbf{v} \cdot \mathbf{u} = 1$.

THEOREM 4. *Given $k$, suppose $\omega_1, \omega_2, \ldots$ are i.i.d. $\sim P_k$. Suppose for each $s \in \mathcal{V}$,*

$$(2.15) \quad \sum_{\mathbf{n} \in \mathcal{V}^*} p_s(\mathbf{n})|\mathbf{n}|^4 < \infty.$$

*Then for any sequence $\lambda_1, \lambda_2, \ldots \in (0,1)$ with*

$$(2.16) \quad \sum_{N=1}^{\infty} \frac{1}{N^2(1-\lambda_N)^2} < \infty,$$

*there is*

$$(2.17) \quad \lim_{N \to \infty} \frac{1-\lambda_N}{N} \sum_{n=1}^{N} S(\omega_n, \lambda_N) = u_k, \qquad P_k\text{-a.s.}$$

REMARK. From Lemma 8, it is seen that,

$$\mathrm{Var}\left(\frac{1-\lambda_N}{N} \sum_{n=1}^{N} S(\omega_n, \lambda_N)\right) = O((1-\lambda_N)^{-1} N^{-1}), \qquad N \to 0.$$

Therefore, if (2.15) is relaxed to $\sum_{\mathbf{n} \in \mathcal{V}^*} p_j(\mathbf{n})|\mathbf{n}|^2 < \infty$ and (2.16) to $(1 - \lambda_N)N \to \infty$, then

$$\frac{1-\lambda_N}{N} \sum_{n=1}^{N} S(\omega_n, \lambda_N) \xrightarrow{P} u_k.$$

**3. Preliminaries.** This section collects some standard results for later use.

LEMMA 1. (a) *(Abel's theorem, cf. [9], Theorems 1.1 and 3.7). Suppose* $\mathbf{w} = (w_1, \ldots, w_V) \in \mathbb{C}^V$ *with* $w_s \neq 0$ *for each* $s$. *If the series*

$$f(\mathbf{z}) := \sum_{\mathbf{n} \in \mathcal{V}^*} a_{\mathbf{n}} z_1^{n_1} \cdots z_V^{n_V}$$

*converges at* $\mathbf{w}$, *then it converges uniformly on any compact subset of* $C_{\mathbf{w}} = \{\mathbf{z} = (z_1, \ldots, z_V) : |z_s| < |w_s|, s \in \mathcal{V}\}$. *The function* $f$ *is analytic on* $C_{\mathbf{w}}$ *and*

$$\frac{\partial f(\mathbf{z})}{\partial z_k} = \sum_{\mathbf{n} \in \mathcal{V}^*} a_{\mathbf{n}} n_k z_k^{n_k-1} \prod_{s \in \mathcal{V} \setminus \{k\}} z_s^{n_s}, \qquad k \in \mathcal{V}.$$

(b) *(cf. [18], Theorem 5.19). Suppose* $f : \Omega \to \mathbb{C}$ *is differentiable on some convex open domain* $\Omega \subset \mathbb{C}^V$ *and continuous on its closure* $\overline{\Omega}$. *Then for any* $\mathbf{z}_1, \mathbf{z}_2 \in \overline{\Omega}$, *there is* $t \in (0, 1)$, *such that*

$$|f(\mathbf{z}_2) - f(\mathbf{z}_1)| \leq |\mathbf{z}_2 - \mathbf{z}_1||\nabla f(t\mathbf{z}_1 + (1-t)\mathbf{z}_2)|.$$

LEMMA 2. *Given* $K \in \mathbb{R}$ *and* $\mathbf{c} = (c_1, \ldots, c_V) \in \mathbb{R}^V$,

$$(3.1) \quad \lim_{t \to 0} \frac{1}{t^2} \sum_{s \in \mathcal{V}} v_s(e^{i(\mathbf{c} \cdot \mathbf{v} - Kt - c_s)t} - 1) = -iK + \frac{1}{2}(\mathbf{c} \cdot \mathbf{v})^2 - \frac{1}{2} \sum_{s \in \mathcal{V}} v_s c_s^2.$$

PROOF. By Taylor's expansion, as $t \to 0$,

$$\sum_{s \in \mathcal{V}} v_s(e^{i(\mathbf{c} \cdot \mathbf{v} - Kt - c_s)t} - 1)$$

$$= \sum_{s \in \mathcal{V}} v_s(i(\mathbf{c} \cdot \mathbf{v} - Kt - c_s)t - \tfrac{1}{2}(\mathbf{c} \cdot \mathbf{v} - Kt - c_s)^2 t^2) + O(t^3)$$

$$= i \sum_{s \in \mathcal{V}} v_s(\mathbf{c} \cdot \mathbf{v} - c_s)t - iK \sum_{s \in \mathcal{V}} v_s t^2 - \tfrac{1}{2} \sum_{s \in \mathcal{V}} v_s(\mathbf{c} \cdot \mathbf{v} - c_s)^2 t^2 + O(t^3).$$

By $\sum_{s \in \mathcal{V}} v_s = 1$, the coefficient of $t$ is 0. Therefore

$$\lim_{t \to 0} \frac{1}{t^2} \sum_{s \in \mathcal{V}} v_s(e^{i(\mathbf{c} \cdot \mathbf{v} - Kt - c_s)t} - 1) = -iK - \frac{1}{2} \sum_{s \in \mathcal{V}} v_s(\mathbf{c} \cdot \mathbf{v} - c_s)^2$$

$$= -iK - \frac{1}{2} \sum_{s \in \mathcal{V}} (v_s(\mathbf{c} \cdot \mathbf{v})^2 - 2 v_s c_s \mathbf{c} \cdot \mathbf{v} + v_s c_s^2),$$

which completes the proof. □

LEMMA 3. *Let $\mathbf{X}_1, \mathbf{X}_2, \ldots$ be random vectors in $\mathbb{R}^V$. Suppose there is a subset $\Omega \subset \mathbb{R}^V$ with Lebesgue measure 0, such that for any $\mathbf{c} \notin \Omega$, $\mathbf{c} \cdot \mathbf{X}_n$ converges in distribution. Then $\mathbf{X}_n \xrightarrow{\mathcal{D}} \mathbf{X}$ for some random vector $\mathbf{X}$ with characteristic function*

$$Ee^{i\mathbf{c}\cdot\mathbf{X}} = \phi(\mathbf{c}) := \lim_{n\to\infty} Ee^{i\mathbf{c}\cdot\mathbf{X}_n} \qquad \forall \mathbf{c} \in \{t\mathbf{x} : t \in \mathbb{R}, \mathbf{x} \notin \Omega\}.$$

PROOF. Because $\mathbb{R}^V \setminus \Omega$ is dense, there exist $\mathbf{c}_1, \ldots, \mathbf{c}_V \notin \Omega$ which are linearly independent, such that $\mathbf{c}_k \cdot X_n$ converges in distribution. The map $T: \mathbf{x} \to (\mathbf{c}_1 \cdot \mathbf{x}, \ldots, \mathbf{c}_V \cdot \mathbf{x})$ is a linear invertible transform on $\mathbb{R}^V$. Because $\{\mathbf{c}_k \cdot \mathbf{X}_n\}$ is tight for each $k \in \mathcal{V}$, so is $\{\mathbf{Y}_n\}$ with $\mathbf{Y}_n = T\mathbf{X}_n$. For any linear transform $A$, $\{A\mathbf{Y}_n\}$ is tight. In particular, with $A = T^{-1}$, $\{\mathbf{X}_n\}$ is tight. Then the characteristic functions $\phi_n(\mathbf{c}) := Ee^{i\mathbf{c}\cdot\mathbf{X}_n}$ are equicontinuous. From $\phi_n(\mathbf{c}) \to \phi(\mathbf{c})$ on a dense subset of $\mathbb{R}^V$, it follows that the convergence holds on the entire $\mathbb{R}^V$, and $\phi$ has a unique continuous extension from $\mathbb{R}^V \setminus \Omega$ to $\mathbb{R}^V$. Now by tightness, $\mathbf{X}_n \xrightarrow{\mathcal{D}} \mathbf{X}$ for some random vector $\mathbf{X}$ and clearly the characteristic of $\mathbf{X}$ has to be $\phi$. □

**4. Limit laws for additive functions.** Suppose $g_k$, $k \in \mathcal{V}$, are real-valued functions on $\mathcal{V}^*$. One can define a function $G$ on the branching trees, such that, for any tree $\omega$ rooted with a particle of type $k$,

$$G(\omega) = g_k(\mathbf{n}) + \sum_{s\in\mathcal{V}} \sum_{j=1}^{n_s} G(\omega_{j,s}),$$

where $(k \to \mathbf{n})$ is the branching rule applied at the root, and $\omega_{j,s}$ is the subtree of $\omega$ rooted with the $j$th particle of type $s$ in $\mathbf{n}$. By recursion, it is easy to check that

$$(4.1) \qquad G(\omega) = \sum_{s\in\mathcal{V}} \sum_{\mathbf{n}\in\mathcal{V}^*} g_s(\mathbf{n}) f(s \to \mathbf{n}; \omega).$$

We will refer to functions of the form (4.1) as "additive" functions for which there is the following theorem.

THEOREM 5. *Assume the same conditions as in Theorem 1. Suppose $G$ is an additive function with $g_1, \ldots, g_V$ satisfying*

$$(4.2) \quad \sum_{s\in\mathcal{V}} v_s \sum_{\mathbf{n}\in\mathcal{V}^*} p_s(\mathbf{n})|g_s(\mathbf{n})| < \infty, \qquad C_g := \sum_{s\in\mathcal{V}} v_s \sum_{\mathbf{n}\in\mathcal{V}^*} p_s(\mathbf{n})g_s(\mathbf{n}) \neq 0.$$

*Given $k \in \mathcal{V}$, let $\phi_k(t) = E_k[e^{itG(\omega)}]$. Then*

$$(4.3) \quad \lim_{t\to 0+} \frac{\phi_k(t)-1}{\sqrt{t}} = L_k := -\frac{u_k[1-i\operatorname{sign}(C_g)]}{\sqrt{H(\mathbf{u})}}\sqrt{|C_g|}, \qquad k \in \mathcal{V}.$$

First, we show that Theorem 1 is implied by the above result.

PROOF OF THEOREM 1. Given $\mathbf{c} = (c_1, \ldots, c_V) \in \mathbb{R}^V$, define $g_s(\mathbf{n}) \equiv c_s$ for any $s \in \mathcal{V}$. Then

$$C_g = \sum_{s \in \mathcal{V}} v_s \sum_{\mathbf{n} \in \mathcal{V}^*} c_s p_s(\mathbf{n}) = \mathbf{c} \cdot \mathbf{v}$$

and

$$G(\omega) = \sum_{s \in \mathcal{V}} \sum_{\mathbf{n} \in \mathcal{V}^*} c_s f(s \to \mathbf{n}; \omega) = \mathbf{c} \cdot \mathbf{f}(\omega).$$

The linear subspace $\mathbf{v}^\perp := \{\mathbf{c} \in \mathbb{R}^V : \mathbf{c} \cdot \mathbf{v} = 0\}$ does not contain $\mathbf{1}$, and hence its Lebesgue measure is 0. By Lemma 3, it is enough to show Theorem 1 for $\mathbf{c} \notin \mathbf{v}^\perp$. For any such $\mathbf{c}$, (4.2) holds, and thus

$$(4.4) \qquad \lim_{t \to 0+} \frac{\phi_k(t) - 1}{\sqrt{t}} = -\frac{u_k[1 - i\,\mathrm{sign}(\mathbf{c} \cdot \mathbf{v})]}{\sqrt{H(\mathbf{u})}} \sqrt{|\mathbf{c} \cdot \mathbf{v}|}, \qquad k \in \mathcal{V}.$$

For $\omega_1, \omega_2, \ldots$ are i.i.d. $\sim P_k$,

$$E\left[\exp\left\{\frac{i}{N^2} \sum_{n=1}^N \mathbf{c} \cdot \mathbf{f}(\omega_n)\right\}\right] = \left[\phi_k\left(\frac{1}{N^2}\right)\right]^N$$

$$= \left(1 + \phi_k\left(\frac{1}{N^2}\right) - 1\right)^N$$

$$= \left(1 + \frac{L_k}{N} + o\left(\frac{1}{N}\right)\right)^N.$$

Then by (4.4), letting $N \to \infty$ leads to

$$\lim_{N \to \infty} E\left[\exp\left\{i\mathbf{c} \cdot \left(\frac{1}{N^2} \sum_{n=1}^N \mathbf{f}(\omega_n)\right)\right\}\right] = e^{L_k} = E[e^{iu_k^2 \mathbf{c} \cdot \mathbf{v} \xi / H(\mathbf{u})}],$$

which completes the proof. □

Following the proof of Theorem 1, we have the following corollary to Theorem 5.

COROLLARY 1. *Suppose $g_s$ satisfies (4.2) and $\omega_1, \omega_2, \ldots \sim P_k$. Then, with $\xi$ being the same as in (2.3),*

$$\frac{1}{N^2} \sum_{n=1}^N G(\omega_n) \xrightarrow{\mathcal{D}} \frac{u_k^2}{H(\mathbf{u})} C_g \xi, \qquad N \to \infty.$$

The rest of the section is devoted to the proof of Theorem 5. First note that, by (2.1) and **u**, **v** being positive,

$$\sum_{\mathbf{n}\in\mathcal{V}^*} p_k(\mathbf{n})(\mathbf{n}\cdot\mathbf{n}) < \infty, \qquad k \in \mathcal{V}. \tag{4.5}$$

LEMMA 4. *Fix $k \in \mathcal{V}$ and function $\theta: \mathcal{V}^* \to \mathbb{C}$, with $|\theta(\mathbf{n})| \leq 1$ for any $\mathbf{n} \in \mathcal{V}^*$. Then $h_k(\mathbf{z}, \theta)$ given below is a well-defined second-order homogeneous polynomial in $\mathbf{z} \in \mathbb{C}^V$:*

$$h_k(\mathbf{z}; \theta) = \sum_{\mathbf{n}\in\mathcal{V}^*} \theta(\mathbf{n}) p_k(\mathbf{n}) \left[\sum_{s\in\mathcal{V}} n_s(n_s-1)z_s^2 + 2\sum_{s<r} n_s n_r z_s z_r\right]. \tag{4.6}$$

*In particular, if $\theta(\mathbf{n}) \equiv 1$, then*

$$h_k(\mathbf{z}; \theta) = h_k(\mathbf{z}) := \sum_{\mathbf{n}\in\mathcal{V}^*} p_k(\mathbf{n})(\mathbf{n}\cdot\mathbf{z})^2 - \sum_{s\in\mathcal{V}} M_{ks} z_s^2. \tag{4.7}$$

*Also define*

$$r_k(\mathbf{z}) = \sum_{\mathbf{n}\in\mathcal{V}^*} \theta(\mathbf{n}) p_k(\mathbf{n}) \prod_{s\in\mathcal{V}} (1+z_s)^{n_s} - q_k(\mathbf{z}), \tag{4.8}$$

*with*

$$q_k(\mathbf{z}) = \sum_{\mathbf{n}\in\mathcal{V}^*} \theta(\mathbf{n}) p_k(\mathbf{n})(1+\mathbf{n}\cdot\mathbf{z}) + \tfrac{1}{2} h_k(\mathbf{z}).$$

*Then $r_k$ is analytic in the interior of $C = \{\mathbf{z} \in \mathbb{C}^V : |1+z_s| < 1, \ s \in \mathcal{V}\}$. Furthermore,*

$$r_k(\mathbf{z}) = o(|\mathbf{z}|^2), \qquad \mathbf{z} \to \mathbf{0}, \ \mathbf{z} \in C. \tag{4.9}$$

PROOF. By (4.5), the summations over $\mathcal{V}^*$ in (4.6) converge, and hence $h_k$ is well defined. Equation (4.7) follows from direct computation. Consider

$$a_k(\mathbf{z}) = \sum_{\mathbf{n}\in\mathcal{V}^*} \theta(\mathbf{n}) p_k(\mathbf{n}) \prod_{s\in\mathcal{V}} z_s^{n_s}.$$

By Lemma 1(a), $a_k$ is analytic in the open domain $D = \{\mathbf{z}: |z_s| < 1, \ s \in \mathcal{V}\}$. In addition, by Lemma 1(a) and (4.5), it is not hard to check that the second-order derivatives of $a_k$ are bounded on $D$, and continuously extend to $\overline{D}$. Then by computation,

$$a_k(\mathbf{1}) = q_k(\mathbf{0}), \qquad \frac{\partial a_k(\mathbf{1})}{\partial z_s} = \frac{\partial q_k(\mathbf{0})}{\partial z_s}, \qquad \frac{\partial^2 a_k(\mathbf{1})}{\partial z_r\, \partial z_s} = \frac{\partial^2 q_k(\mathbf{0})}{\partial z_r\, \partial z_s}, \qquad r,s \in \mathcal{V}.$$

Consequently, the first- and second-order derivatives of

$$r_k(\mathbf{z}) = a_k(\mathbf{z}+\mathbf{1}) - \sum_{\mathbf{n}\in\mathcal{V}^*} \theta(\mathbf{n}) p_k(\mathbf{n})(1+\mathbf{n}\cdot\mathbf{z}) - \tfrac{1}{2} h_k(\mathbf{z})$$

are bounded on $C$ and have continuous extension to $\overline{C}$, such that

$$\lim_{\mathbf{z}\to 0} r_k(\mathbf{z}) = \lim_{\mathbf{z}\to 0} \frac{\partial r_k(\mathbf{z})}{\partial z_s} = \lim_{\mathbf{z}\to 0} \frac{\partial^2 r_k(t\mathbf{z})}{\partial z_r \partial z_s} = 0.$$

Apply Lemma 1(b) twice, once to $r_k$ and once to its derivatives. Then it can be seen that, for any $\mathbf{z} \in \overline{C}$, there is $t \in (0,1)$, such that

$$|r_k(\mathbf{z})| = |r_k(\mathbf{z}) - r_k(\mathbf{0})| \leq |\mathbf{z}|^2 \max_{r,s \in \mathcal{V}} \left| \frac{\partial^2 r_k(t\mathbf{z})}{\partial z_r \partial z_s} \right|,$$

which leads to (4.9). $\square$

LEMMA 5. *Recall function $H$ defined in* (2.4). *For the function, we have*

(4.10)
$$H(\mathbf{z}) = \sum_{s \in \mathcal{V}} v_s h_s(\mathbf{z}) = E_Q(\mathbf{X} \cdot \mathbf{z})^2 - \sum_{s \in \mathcal{V}} v_s z_s^2,$$

$H(\mathbf{z}) \not\equiv 0$ *and all its coefficients are nonnegative.*

PROOF. Because $\mathbf{v} = \mathbf{v}M$, from (2.4), it is not hard to see that the equalities in (4.10) hold. For each $k \in \mathcal{V}$, all the coefficients of $h_k$ are nonnegative and $v_k > 0$. Therefore, all the coefficients of $H$ are nonnegative as well. If $H(\mathbf{z}) \equiv 0$, then there must be $h_k \equiv 0$, $k \in \mathcal{V}$. From (4.7), this implies that $p_k(\mathbf{n}) > 0$ only if $\mathbf{n} = \mathbf{0}$ or $\mathbf{n} = \mathbf{e}_s = (\varepsilon_{s1}, \ldots, \varepsilon_{sV})$, for some $s \in \mathcal{V}$, with $\varepsilon_{sr} = 0$ if $r \neq s$ and 1 otherwise. Therefore $M_{rs} = p_r(\mathbf{e}_s)$ and $\sum_{s \in \mathcal{V}} p_r(\mathbf{e}_s) u_s = u_r$. By choosing $r$ such that $u_r = \max\{u_s, s \in \mathcal{V}\}$, it is seen that $u_s = u_r$ for all $s \in \mathcal{V}$ and $\sum_{s \in \mathcal{V}} p_r(\mathbf{e}_s) = 1$. Thus, almost surely, each particle produces exactly one offspring, leading to a nonterminating process, which is a contradiction. $\square$

PROOF OF THEOREM 5. The following recursive relations hold:

(4.11) $$\phi_k(t) = \sum_{\mathbf{n} \in \mathcal{V}^*} e^{itg_k(\mathbf{n})} p_k(\mathbf{n}) \prod_{\in \mathcal{V}} [\phi_s(t)]^{n_s}, \qquad k \in \mathcal{V}.$$

Let $\Delta_k(t) = \phi_k(t) - 1$, and $\mathbf{\Delta}_t = (\Delta_1(t), \ldots, \Delta_V(t))$. By (4.11),

$$\Delta_k(t) = -1 + \sum_{\mathbf{n} \in \mathcal{V}^*} e^{itg_k(\mathbf{n})} p_k(\mathbf{n}) \prod_{s \in \mathcal{V}} [1 + \Delta_s(t)]^{n_s}$$

$$= -1 + \sum_{\mathbf{n} \in \mathcal{V}^*} e^{itg_k(\mathbf{n})} p_k(\mathbf{n})(1 + \mathbf{n} \cdot \mathbf{\Delta}_t) + \tfrac{1}{2} h_k(\mathbf{\Delta}_t; e^{itg_k}) + r_k(\mathbf{\Delta}_t)$$

$$= \sum_{\mathbf{n} \in \mathcal{V}^*} (e^{itg_k(\mathbf{n})} - 1) p_k(\mathbf{n}) + \sum_{\mathbf{n} \in \mathcal{V}^*} e^{itg_k(\mathbf{n})} p_k(\mathbf{n}) \mathbf{n} \cdot \mathbf{\Delta}_t$$

(4.12)
$$\quad + \tfrac{1}{2} h_k(\mathbf{\Delta}_t; e^{itg_k}) + r_k(\mathbf{\Delta}_t)$$

$$= \sum_{\mathbf{n} \in \mathcal{V}^*} (e^{itg_k(\mathbf{n})} - 1) p_k(\mathbf{n})(1 + \mathbf{n} \cdot \mathbf{\Delta}_t) + \sum_{s \in \mathcal{V}} M_{ks} \Delta_s(t)$$

$$\quad + \tfrac{1}{2} h_k(\mathbf{\Delta}_t; e^{itg_k}) + r_k(\mathbf{\Delta}_t),$$

where $h_k$ is defined as in (4.6), with function $\theta = e^{itg_k} : \mathbf{n} \to e^{itg_k(\mathbf{n})}$, and $r_k$ is defined as in (4.8).

We will use (4.12) to prove (4.3). It is enough to show that, for any $t_n \to 0+$, there is a subsequence $\tau_n$ of $t_n$, such that $\tau_n^{-1/2} \mathbf{\Delta}_{\tau_n}$ converge to $\mathbf{L} = (L_1, \ldots, L_V)$, with $L_k$ given in (4.3).

Fix an arbitrary $t_n \to 0+$. We first show that, when $n$ is large enough, $\mathbf{\Delta}_{t_n} \neq 0$. Indeed, if this is not the case, then (4.12) implies

$$\sum_{\mathbf{n} \in \mathcal{V}^*} (e^{itg_k(\mathbf{n})} - 1) p_k(\mathbf{n}) = 0 \qquad \text{for } t = t_n, \ n \text{ large enough.}$$

Then by dominated convergence, for all $k \in \mathcal{V}$,

$$\sum_{\mathbf{n} \in \mathcal{V}^*} g_k(\mathbf{n}) p_k(\mathbf{n}) = \lim_{n \to \infty} \frac{1}{it_n} \sum_{\mathbf{n} \in \mathcal{V}^*} (e^{it_n g_k(\mathbf{n})} - 1) p_k(\mathbf{n}) = 0,$$

contradicting (4.2). Then there is a subsequence $\{t'_n\} \subset \{t_n\}$ as well as a subset $\mathcal{V}' \subset \mathcal{V}$, such that $\Delta_s(t'_n) \neq 0$ for all $s \in \mathcal{V}'$ while $\Delta_s(t'_n) = 0$ for all $s \notin \mathcal{V}'$.

We show that $\exists k_0 \in \mathcal{V}'$ and $\{\tau_n\} \subset \{t'_n\}$, such that

$$(4.13) \qquad \xi_s = \lim_{n \to \infty} \frac{\Delta_s(\tau_n)}{\Delta_{k_0}(\tau_n)}$$

exists for any $s \in \mathcal{V}$. First, for any $s \notin \mathcal{V}'$, (4.13) is clear. Suppose $\mathcal{V}' = \{s_1, \ldots, s_m\}$. If $m = 1$, (4.13) is also obvious. If $m > 1$, then

$$\liminf_{n \to \infty} \left| \frac{\Delta_{s_1}(t_n)}{\Delta_{s_2}(t_n)} \right| < \infty \quad \text{or} \quad \liminf_{n \to \infty} \left| \frac{\Delta_{s_2}(t_n)}{\Delta_{s_1}(t_n)} \right| < \infty.$$

Assuming the first one, for some $\{t''_n\} \subset \{t'_n\}$, $\frac{\Delta_{s_1}(t''_n)}{\Delta_{s_2}(t''_n)}$ converges. By induction, there exists $k_0 \in \mathcal{V}' \setminus \{s_1\}$ as well as $\{\tau_n\} \subset \{t''_n\}$, such that $\frac{\Delta_s(\tau_n)}{\Delta_{k_0}(\tau_n)}$ converges for all $s \in \mathcal{V}' \setminus \{s_1\}$. In particular, $\frac{\Delta_{s_2}(\tau_n)}{\Delta_{k_0}(\tau_n)}$ converges, implying that $\frac{\Delta_{s_1}(\tau_n)}{\Delta_{k_0}(\tau_n)}$ converges as well.

With $k_0$ and $\tau_0$ being fixed such that (4.13) holds, denote $\boldsymbol{\xi} = (\xi_1, \ldots, \xi_V)$. Clearly $\xi_{k_0} = 1$. To get the other $\xi_l$, first consider the asymptotics of $\mathbf{\Delta}_t$. Rewrite (4.12) to get

$$(4.14) \quad \begin{aligned} \Delta_k(t) - \sum_{s \in \mathcal{V}} M_{ks} \Delta_s(t) + \sum_{\mathbf{n} \in \mathcal{V}^*} (1 - e^{itg_k(\mathbf{n})}) p_k(\mathbf{n})(1 + \mathbf{n} \cdot \mathbf{\Delta}_t) \\ = \tfrac{1}{2} h_k(\mathbf{\Delta}_t; e^{itg_k}) + r_k(\mathbf{\Delta}_t). \end{aligned}$$

Multiply both sides of (4.14) by $v_k$ and then sum over $\mathcal{V}$. Because $\mathbf{v} = \mathbf{v}M$,

$$\sum_{s \in \mathcal{V}} v_s \sum_{\mathbf{n} \in \mathcal{V}^*} (1 - e^{itg_s(\mathbf{n})}) p_s(\mathbf{n})(1 + \mathbf{n} \cdot \mathbf{\Delta}_t)$$

(4.15)
$$= \tfrac{1}{2}\sum_{s\in\mathcal{V}} v_s h_s(\boldsymbol{\Delta}_t; e^{itg_s}) + \sum_{s\in\mathcal{V}} v_s r_s(\boldsymbol{\Delta}_t).$$

Divide both sides of (4.15) by $\Delta_{k_0}^2(t)$ and let $t\to 0$ through $\tau_n$. Because each $h_k(\mathbf{z}, e^{itg_s})$ is a second-order homogeneous polynomial, by dominated convergence,
$$\lim_{n\to\infty} \frac{1}{\Delta_{k_0}^2(\tau_n)} h_s(\boldsymbol{\Delta}_{\tau_n}; e^{i\tau_n g_s}) = h_s(\boldsymbol{\xi}).$$

Since $|\phi_s(t)| = |1+\Delta_s(t)| \leq 1$ for all $s\in\mathcal{V}$, and $\boldsymbol{\Delta}_t \to 0$ as $t\to 0$, (4.9) leads to

(4.16)
$$r_s(\boldsymbol{\Delta}_t) = o(|\boldsymbol{\Delta}_t|^2), \qquad t\to 0.$$

On the other hand, by (4.2) and dominated convergence,
$$\lim_{n\to\infty} \frac{1}{\tau_n} \sum_{s\in\mathcal{V}} v_s \sum_{\mathbf{n}\in\mathcal{V}^*} (1 - e^{i\tau_n g_s(\mathbf{n})}) p_s(\mathbf{n})(1+\mathbf{n}\cdot\boldsymbol{\Delta}_t)$$
$$= -i \sum_{s\in\mathcal{V}} \sum_{\mathbf{n}\in\mathcal{V}^*} v_s p_s(\mathbf{n}) g_s(\mathbf{n}) = -iC_g.$$

Combining the above results,

(4.17)
$$\lim_{n\to\infty} \frac{\tau_n}{\Delta_{k_0}^2(\tau_n)} = \frac{1}{-2iC_g}\sum_{s\in\mathcal{V}} v_s h_s(\boldsymbol{\xi}) = \frac{1}{-2iC_g} H(\boldsymbol{\xi})$$

and hence

(4.18)
$$\lim_{n\to\infty} \frac{\tau_n}{\Delta_{k_0}(\tau_n)} = 0.$$

Divide both sides of (4.14) by $\Delta_{k_0}(t)$, and then let $t\to 0$ along $\tau_n$. By (4.16)–(4.18) and $h_k(\boldsymbol{\Delta}_t) = o(\Delta_{k_0}(t))$, there is $\boldsymbol{\xi} = M\boldsymbol{\xi}$. Therefore, $\boldsymbol{\xi}$ is an eigenvector corresponding to the simple eigenvalue 1 of $M$, and thus $\boldsymbol{\xi}$ is some constant times $\mathbf{u}$. Since $\xi_{k_0} = 1$, by comparing with (1.1), we get

(4.19)
$$\xi_s = \lim_{n\to\infty} \frac{\Delta_s(\tau_n)}{\Delta_{k_0}(\tau_n)} = \frac{u_s}{u_{k_0}} > 0, \qquad s\in\mathcal{V}, \qquad \boldsymbol{\xi} = \frac{\mathbf{u}}{u_{k_0}},$$

which, together with (4.10) and (4.18), leads to
$$\lim_{n\to\infty} \frac{\Delta_{k_0}^2(\tau_n)}{\tau_n} = -\frac{2C_g u_{k_0}^2 i}{H(\mathbf{u})},$$

where $H(\mathbf{u}) > 0$ is because $\mathbf{u}$ is positive and $H \not\equiv 0$ with all its coefficients nonnegative. From the limit it is seen that, as $n\to\infty$, $\Delta_{k_0}^2(\tau_n)/\tau_n$ at most has two cluster points. Because the real part of $\Delta_{k_0}(t) = E_{k_0}[e^{itG(\omega)}] - 1$ is nonpositive,

$$\lim_{n \to \infty} \frac{\Delta_{k_0}(\tau_n)}{\sqrt{\tau_n}} = -\frac{u_{k_0}(1 - i\,\text{sign}(C_g))}{\sqrt{H(\mathbf{u})}}\sqrt{|C_g|}. \tag{4.20}$$

Combining the above limit and (4.19), we get that, for any $s \in \mathcal{V}$,

$$\lim_{n \to \infty} \frac{\Delta_s(\tau_n)}{\sqrt{\tau_n}} = -\frac{u_s(1 - i\,\text{sign}(C_g))}{\sqrt{H(\mathbf{u})}}\sqrt{|C_g|}.$$

This completes the proof. □

**5. Limit laws for the relative frequencies of types.** This section is devoted to the proof of Theorem 2. By Lemma 3, it is enough to establish the result for $K \neq 0$. Given $\mathbf{c} = (c_1, \ldots, c_V) \in \mathbb{R}^V$ and $K \in \mathbb{R} \setminus \{0\}$,

$$\frac{1}{N}\sum_{n=1}^{N}(\mathbf{c} \cdot \mathbf{f}(\omega_n) - (\mathbf{c} \cdot \mathbf{v})|\omega_n|) + \frac{K}{N^2}\sum_{n=1}^{N}|\omega_n|$$

$$= \frac{1}{N}\sum_{n=1}^{N}\mathbf{c} \cdot \mathbf{f}(\omega_n) - \frac{1}{N}\sum_{n=1}^{N}\left(\mathbf{c} \cdot \mathbf{v} - \frac{K}{N}\right)|\omega_n|.$$

Note that, in Theorem 2, it is assumed that $\omega_1, \omega_2, \ldots$ are i.i.d. $\sim P_k$. Therefore, letting

$$\theta_k(t) = E_k[e^{it(\mathbf{c} \cdot \mathbf{f}(\omega) - C_t|\omega|)}], \qquad C_t = \mathbf{c} \cdot \mathbf{v} - Kt, \ k \in \mathcal{V} \tag{5.1}$$

we need to find the limit

$$\lim_{N \to \infty}\left[\theta_k\left(\frac{1}{N}\right)\right]^N = \lim_{t \to 0+}[\theta_k(t)]^{1/t}.$$

Following the previous section, let

$$\Delta_k(t) = \theta_k(t) - 1, \qquad \mathbf{\Delta}_t = (\Delta_1(t), \ldots, \Delta_V(t)).$$

Then, as in the proof of Theorem 1, provided $L_k := \lim_{t \to 0+}\frac{1}{t}\Delta_k(t)$ exists and

$$L_k = z(\mathbf{c}, K)u_k + i\eta_k, \qquad k \in \mathcal{V},$$

$\lim_{t \to 0+}[\theta_k(t)]^{1/t} = e^{L_k}$. To this end, it is enough to show that, for any $t_n \to 0+$, there is a subsequence $\{\tau_n\} \subset \{t_n\}$, such that $\tau_n^{-1}\mathbf{\Delta}_{\tau_n} \to \mathbf{L} = (L_1, \ldots, L_V) \in \mathbb{C}^V$.

The functions $\theta_k(t)$, $k \in \mathcal{V}$, have the following recursive relations:

$$\theta_k(t) = \sum_{\mathbf{n} \in \mathcal{V}^*} p_k(\mathbf{n})e^{i(c_k - C_t)t}\prod_{s \in \mathcal{V}}(\theta_s(t))^{n_s}.$$

Therefore,
$$e^{i(C_t-c_k)t}(1+\Delta_k(t)) = \sum_{\mathbf{n}\in\mathcal{V}^*} p_k(\mathbf{n})\prod_{s\in\mathcal{V}}(1+\Delta_s(t))^{n_s}.$$

Then similarly to (4.12), it can be shown that

$$(5.2) \Delta_k(t) + (e^{i(C_t-c_k)t}-1)(1+\Delta_k(t)) = \sum_{s\in\mathcal{V}} M_{ks}\Delta_s(t) + \tfrac{1}{2}h_k(\boldsymbol{\Delta}_t) + r_k(\boldsymbol{\Delta}_t),$$

where $h_k$ and $r_k$ are defined as in (4.7) and (4.8), respectively. Multiply both sides of (5.2) by $v_k$ and take the sum over $k \in \mathcal{V}$. By $v_1 + \cdots + v_V = 1$, $\mathbf{v} = \mathbf{v}M$ and (4.10), one gets

$$(5.3) \sum_{s\in\mathcal{V}} v_s(e^{i(C_t-c_s)t}-1) + \sum_{s\in\mathcal{V}} v_s(e^{i(C_t-c_s)t}-1)\Delta_s(t) = \tfrac{1}{2}H(\boldsymbol{\Delta}_t) + r(\boldsymbol{\Delta}_t),$$

where $r(\mathbf{z}) := \sum_{k\in\mathcal{V}} v_k r_k(\mathbf{z})$. Fix $t_n \to 0+$. Following (4.13), we want to find $k_0 \in \mathcal{V}$ and $\{t'_n\} \subset \{t_n\}$ such that, for any $s \in \mathcal{V}$,

$$(5.4) \qquad \xi_s = \lim_{n\to\infty} \frac{\Delta_s(t'_n)}{\Delta_{k_0}(t'_n)}$$

exists. First, we show that, when $t > 0$ is small enough, $\boldsymbol{\Delta}_t \neq 0$. Suppose this is not the case. Then (5.3) implies that, for a sequence of $t \to 0$, $\sum_{s\in\mathcal{V}} v_s(e^{i(C_t-c_s)t}-1) = 0$. According to Lemma 2, this leads to

$$\lim_{t\to 0} \frac{1}{t^2} \sum_{s\in\mathcal{V}} v_s(e^{i(C_t-c_s)t}-1) = -2iK + (\mathbf{c}\cdot\mathbf{v})^2 - \sum_{s\in\mathcal{V}} v_s c_s^2 = 0.$$

In particular, $K = 0$, which is a contradiction. Thus $\boldsymbol{\Delta}_t \neq 0$ for $t > 0$ small. Then we can find a subsequence $\{t''_n\} \subset \{t_n\}$ and a subset $\mathcal{V}' \subset \mathcal{V}$, such that $\Delta_s(t''_n) \neq 0$ for all $s \in \mathcal{V}'$ while $\Delta_s(t''_n) = 0$ for all $s \notin \mathcal{V}'$. With an argument similar to the one for (4.13), there are $\{t'_n\} \subset \{t''_n\}$ and $k_0 \in \mathcal{V}'$ such that (5.4) holds for $s \in \mathcal{V}'$. For $s \notin \mathcal{V}'$, (5.4) clearly holds with $\xi_s = 0$. Thus (5.4) holds for all $s \in \mathcal{V}$.

LEMMA 6. *With $k_0$ chosen as above,*

$$\limsup_{n\to\infty} \left|\frac{\Delta_{k_0}(t'_n)}{t'_n}\right| < \infty.$$

Assume for the moment that Lemma 6 is true. Fix $\mathbf{c} \in \mathbb{R}^V$ with $c_1, \ldots, c_V$ not all equal. By the lemma, there is $\{\tau_n\} \subset \{t_n\}$, such that $\tau_n^{-1}\Delta_{k_0}(\tau_n)$ converge. Therefore, by (5.4), there is $\boldsymbol{\lambda} \in \mathbb{C}^V$, such that $\lim_{n\to\infty} \boldsymbol{\Delta}_{\tau_n}/\tau_n = \boldsymbol{\lambda}$. Divide (5.2) by $t$ and then let $t \to 0$ through $\tau_n$ to get

$$(5.5) \; \boldsymbol{\lambda}^t + i((\mathbf{c}\cdot\mathbf{v})\mathbf{1}^t - \mathbf{c}^t) = M\boldsymbol{\lambda}^t \quad \Longrightarrow \quad i(\mathbf{1}^t\mathbf{v} - \mathbf{v}^t)\mathbf{c}^t = (M-I)\boldsymbol{\lambda}^t.$$

By (2.7), if $\boldsymbol{\eta}^t = \Lambda \mathbf{c}^t$, then $\boldsymbol{\lambda} := i\boldsymbol{\eta}$ is a solution to (5.5). Moreover, since 1 is a simple eigenvalue of $M$, any solution to (5.5) can be written as $z\mathbf{u} + i\boldsymbol{\eta}$, $z \in \mathbb{C}$. As a result, $\lim \tau_n^{-1} \boldsymbol{\Delta}_{\tau_n} = z\mathbf{u} + i\boldsymbol{\eta}$ for some $z \in \mathbb{C}$ which is to be found.

Divide both sides of (5.3) by $t^2$ and let $t \to 0$ through $\tau_n$. Since $H$ is a homogeneous polynomial of order 2 and, according to (4.9), $R(\boldsymbol{\Delta}_t) = o(|\boldsymbol{\Delta}_t|^2)$, then (3.1) implies

$$-iK + \tfrac{1}{2}(\mathbf{c}\cdot\mathbf{v})^2 - \tfrac{1}{2}\sum_{s\in\mathcal{V}} v_s c_s^2 + i\sum_{s\in\mathcal{V}} v_s(\mathbf{c}\cdot\mathbf{v} - c_s)(zu_s + i\eta_s) = \tfrac{1}{2}H(z\mathbf{u} + i\boldsymbol{\eta}).$$

By some calculation, it can be seen that $z$ is a solution to the equation

$$H(\mathbf{u})z^2 + 2Azi + B + 2Ki = 0,$$

where

$$A = E_Q[(\mathbf{X}\cdot\mathbf{u})(\mathbf{X}\cdot\boldsymbol{\eta})] - \sum_{s\in\mathcal{V}} v_s u_s(\eta_s - c_s) - (\mathbf{v}\cdot\mathbf{c})(\mathbf{v}\cdot\mathbf{u})$$

$$= \mathrm{Cov}_Q(\mathbf{X}\cdot\mathbf{u}, \mathbf{X}\cdot\boldsymbol{\eta}) + E_Q(\mathbf{X}\cdot\mathbf{u})E_Q(\mathbf{X}\cdot\boldsymbol{\eta}) - \sum_{s\in\mathcal{V}} v_s u_s(\eta_s - c_s) - \mathbf{v}\cdot\mathbf{c},$$

$$B = -H(\boldsymbol{\eta}) - 2\sum_{s\in\mathcal{V}} v_s c_s \eta_s + 2(\mathbf{v}\cdot\mathbf{c})(\mathbf{v}\cdot\boldsymbol{\eta}) - (\mathbf{c}\cdot\mathbf{v})^2 + \sum_{s\in\mathcal{V}} v_s c_s^2$$

$$= -E_Q(\mathbf{X}\cdot\boldsymbol{\eta})^2 + \sum_{s\in\mathcal{V}} v_s \eta_s^2 - 2\sum_{s\in\mathcal{V}} v_s c_s \eta_s + 2(\mathbf{v}\cdot\mathbf{c})(\mathbf{v}\cdot\boldsymbol{\eta}) - (\mathbf{c}\cdot\mathbf{v})^2 + \sum_{s\in\mathcal{V}} v_s c_s^2$$

$$= -\mathrm{Var}_Q(\mathbf{X}\cdot\boldsymbol{\eta}) - (E_Q(\mathbf{X}\cdot\boldsymbol{\eta}))^2 + 2(\mathbf{v}\cdot\mathbf{c})(\mathbf{v}\cdot\boldsymbol{\eta}) + \sum_{s\in\mathcal{V}} v_s(\eta_s - c_s)^2 - (\mathbf{c}\cdot\mathbf{v})^2.$$

On the one hand,

$$E_Q(\mathbf{X}) = \left(\sum_{k\in\mathcal{V}}\sum_{\mathbf{n}\in\mathcal{V}^*} v_k p_k(\mathbf{n}) n_1, \ldots, \sum_{k\in\mathcal{V}}\sum_{\mathbf{n}\in\mathcal{V}^*} v_k p_k(\mathbf{n}) n_V\right)$$

$$= \left(\sum_{k\in\mathcal{V}} v_k M_{k1}, \ldots, \sum_{k\in\mathcal{V}} v_k M_{kV}\right) = \mathbf{v},$$

and on the other,

$$\mathbf{v}\Lambda = \sum_{n=0}^{\infty} \mathbf{v}(M^n - \mathbf{u}^t\mathbf{v})(I - \mathbf{1}^t\mathbf{v}) = 0,$$

and hence $\mathbf{v}\cdot\boldsymbol{\eta} = \mathbf{v}\Lambda\mathbf{c}^t = 0$. Therefore, $A$ and $B$ can be expressed as in (2.10) and (2.11).

The roots of the equation are

$$\frac{1}{H(\mathbf{u})}\left(-Ai \pm \sqrt{-(A^2 + BH(\mathbf{u})) - 2KH(\mathbf{u})i}\right).$$

Since the real part of $\Delta_k(t) = \phi_k(t) - 1$ is nonpositive, $z$ is the one with nonpositive real part.

PROOF OF LEMMA 6. If the claim is false, then there exists $\{\tau_n\} \subset \{t'_n\}$ such that
$$\lim_{n \to \infty} \frac{\tau_n}{\Delta_{k_0}(\tau_n)} = 0.$$

First, divide both sides of (5.2) by $\Delta_{k_0}(t)$ and let $t \to 0$ through $\tau_n$. By (5.4), we get $\boldsymbol{\xi} = M\boldsymbol{\xi}$. Since $\xi_{k_0} = 1$, $\boldsymbol{\xi} = u_{k_0}^{-1}\mathbf{u}$. In particular, $H(\boldsymbol{\xi}) > 0$. Thus, by dividing both sides of (5.3) by $\Delta_{k_0}(t)^2$, and letting $t \to 0$ through $\tau_n$, it is see that $0 = H(\boldsymbol{\xi})$, which is a contradiction. □

COROLLARY 2. *For A and B, we have $A^2 \leq -H(\mathbf{u})B$, that is,*
$$\left(\mathrm{Cov}_Q(\mathbf{X} \cdot \mathbf{u}, \mathbf{X} \cdot \boldsymbol{\eta}) - \sum_{s \in \mathcal{V}} v_s u_s(\eta_s - c_s) - \mathbf{c} \cdot \mathbf{v}\right)^2$$
$$\leq H(\mathbf{u})\left(\mathrm{Var}_Q(\mathbf{X} \cdot \boldsymbol{\eta}) - \sum_{s \in \mathcal{V}} v_s(c_s - \eta_s)^2 + (\mathbf{c} \cdot \mathbf{v})^2\right).$$

PROOF. Indeed, if $A^2 > -H(\mathbf{u})B$, then $-(A^2 + H(\mathbf{u})B) < 0$. Since $z$ has to be the one of
$$\frac{1}{H(\mathbf{u})}\left(-Ai \pm \sqrt{-(A^2 + BH(\mathbf{u})) - 2KH(\mathbf{u})i}\right)$$
with negative real part, with $\mathbf{c}$ being fixed, $z$ is not continuous at $K = 0$, a contradiction. □

**6. Limit laws for the relative frequencies of branching rules.** In this section, we prove Theorem 3. With a little abuse of notation, denote
$$\mathbf{F}(\omega) = (f(j \to \mathbf{n}; \omega), \mathbf{n} \in \mathcal{V}^*), \qquad \mathbf{q} = (p_j(\mathbf{n}), \mathbf{n} \in \mathcal{V}^*),$$
and for $\mathbf{c}, \mathbf{q} \in \mathbb{R}^{\mathcal{V}^*}$, provided $\sum_{\mathbf{n} \in \mathcal{V}^*} |c_\mathbf{n} q_\mathbf{n}| < \infty$, $\mathbf{c} \cdot \mathbf{q} = \sum_{\mathbf{n} \in \mathcal{V}^*} c_\mathbf{n} q_\mathbf{n}$.

Given $\mathbf{c} \in \mathbb{R}^{\mathcal{V}^*}$, with $c_\mathbf{n} = 0$ for all but a finite number of $\mathbf{n}$ satisfying $p_j(\mathbf{n}) > 0$, and $K \in \mathbb{R}$, consider $E[e^{i\zeta_N(\mathbf{c}, K)}]$, with the random variable
$$\zeta_N(\mathbf{c}, K) = \frac{1}{N}\sum_{n=1}^{N}(\mathbf{c} \cdot \mathbf{F}(\omega_n) - (\mathbf{c} \cdot \mathbf{q})f(j; \omega_n)) + \frac{K}{N^2}\sum_{n=1}^{N} f(j; \omega_n)$$
$$= \frac{1}{N}\sum_{n=1}^{N} \mathbf{c} \cdot \mathbf{F}(\omega_n) - \frac{1}{N}\sum_{n=1}^{N}\left(\mathbf{c} \cdot \mathbf{q} - \frac{K}{N}\right)f(j; \omega_n).$$

Then, to prove Theorem 3, it is enough to show that

$$E[e^{i\zeta_N(\mathbf{c},K)}] \to e^{z(\mathbf{c},K)u_k},$$

such that $z(\mathbf{c}, K)$ is the solution to

(6.1) $$z^2 + \frac{2v_j A i}{H(\mathbf{u})} z + \frac{v_j}{H(\mathbf{u})}(B + 2Ki) = 0,$$

with

$$A = \sum_{\mathbf{n}\in\mathcal{V}^*} c_{\mathbf{n}} p_j(\mathbf{n})(\mathbf{n}\cdot\mathbf{u}) - (\mathbf{c}\cdot\mathbf{q})u_j,$$

$$B = \sum_{\mathbf{n}\in\mathcal{V}^*} p_j(\mathbf{n})c_{\mathbf{n}}^2 - (\mathbf{c}\cdot\mathbf{q})^2.$$

PROOF OF THEOREM 3. As in the proof of Theorem 2, it is enough to show the above limit for $K \neq 0$. Let

(6.2) $$\psi_j(t) = E_j[e^{it(\mathbf{c}\cdot\mathbf{F}(j\to n;\omega) - C_t f(j;\omega))}] \quad \text{with } C_t = \mathbf{c}\cdot\mathbf{q} - Kt.$$

Then $E[e^{i\zeta_N(\mathbf{c},K)}] = [\psi_j(\frac{1}{N})]^{N^2}$. By recursion,

$$e^{itC_t}\psi_j(t) = \sum_{\mathbf{n}\in\mathcal{V}^*} p_j(\mathbf{n})e^{itc_{\mathbf{n}}} \prod_{s\in\mathcal{V}}[\psi_s(t)]^{n_s},$$

$$\psi_k(t) = \sum_{\mathbf{n}\in\mathcal{V}^*} p_k(\mathbf{n}) \prod_{s\in\mathcal{V}}[\psi_s(t)]^{n_s} \quad \text{for } k \neq j.$$

Let $\Delta_s(t) = \psi_s(t) - 1$, $s \in \mathcal{V}$, and $\boldsymbol{\Delta}_t = (\Delta_1(t), \ldots, \Delta_V(t))$. Then

(6.3) $$\begin{aligned}\Delta_j(t) + (e^{itC_t} - 1)\psi_j(t) \\ = \sum_{s\in\mathcal{V}} M_{js}\Delta_s(t) + \tfrac{1}{2}h_j(\boldsymbol{\Delta}_t) + \sum_{\mathbf{n}\in\mathcal{V}^*} p_j(\mathbf{n})(e^{ic_{\mathbf{n}}t} - 1) \\ + r_j(\boldsymbol{\Delta}_t) + \sum_{\mathbf{n}\in\mathcal{V}^*} p_j(\mathbf{n})(e^{ic_{\mathbf{n}}t} - 1)(\mathbf{n}\cdot\boldsymbol{\Delta}_t + T_{\mathbf{n}}(\boldsymbol{\Delta}_t)),\end{aligned}$$

(6.4) $$\Delta_k(t) = \sum_{s\in\mathcal{V}} M_{ks}\Delta_s(t) + \tfrac{1}{2}h_k(\boldsymbol{\Delta}_t) + r_k(\boldsymbol{\Delta}_t) \quad \text{for } k \neq j,$$

where $h_k$, $k \in \mathcal{V}$, are defined as in (4.7), and $r_k$ as in (4.8). By (4.9), $r_k(\boldsymbol{\Delta}_t) = o(|\boldsymbol{\Delta}_t|^2)$ as $t \to 0$. On other hand, $T_{\mathbf{n}}(\mathbf{z}) = \prod_{s\in\mathcal{V}}(1+\mathbf{z})^{n_s} - \mathbf{n}\cdot\mathbf{z}$ is a polynomial of order greater than or equal to 2. Note that in the last sum in (6.3), since $c_{\mathbf{n}} = 0$ for all but a finite number of $\mathbf{n}$, only a finite number of summands is nonzero. Thus, following the steps in Sections 4 and 5, multiply both sides of (6.3) by $v_j$, both sides of (6.4) by $v_k$, add them up and use $\mathbf{v} = \mathbf{v}M$ to get

$$v_j(e^{iC_t t} - 1)\Delta_j(t)$$

$$= \tfrac{1}{2}H(\boldsymbol{\Delta}_t) + v_j\left[\sum_{\mathbf{n}\in\mathcal{V}^*} p_j(\mathbf{n})(e^{ic_\mathbf{n} t} - 1) - (e^{itC_t} - 1)\right] + r(\boldsymbol{\Delta}_t)$$

(6.5)
$$+ v_j \sum_{\mathbf{n}\in\mathcal{V}^*} p_j(\mathbf{n})(e^{ic_\mathbf{n} t} - 1)(\mathbf{n}\cdot\boldsymbol{\Delta}_t + T_\mathbf{n}(\boldsymbol{\Delta}_t))$$

$$= \tfrac{1}{2}H(\boldsymbol{\Delta}_t) - v_j \sum_{\mathbf{n}\in\mathcal{V}^*} p_j(\mathbf{n})e^{ic_\mathbf{n} t}(e^{it(\mathbf{c}\cdot\mathbf{q}-Kt-c_\mathbf{n} t)} - 1)$$

$$+ v_j \sum_{\mathbf{n}\in\mathcal{V}^*} q_\mathbf{n}(e^{ic_\mathbf{n} t} - 1)(\mathbf{n}\cdot\boldsymbol{\Delta}_t) + o(\boldsymbol{\Delta}_t^2),$$

where $r = \sum_{s\in\mathcal{V}} v_s r_s$. Similar to the argument following (5.3), given any $t_n \to 0+$, there exist $k_0 \in \mathcal{V}$ and a subsequence $\{t'_n\} \subset \{t_n\}$, such that $\Delta_{k_0}(t'_n) \ne 0$, and

(6.6)
$$\lim_{n\to\infty} \frac{\boldsymbol{\Delta}_{t'_n}}{\Delta_{k_0}(t'_n)} = \boldsymbol{\xi},$$

for some $\boldsymbol{\xi} \in \mathbb{C}^V$. We need the following bounds. $\square$

LEMMA 7. *With $k_0$ chosen as above,*

(6.7)
$$\limsup_{n\to\infty}\left|\frac{\Delta_{k_0}(t'_n)}{t'_n}\right| < \infty.$$

Assume Lemma 7 is true for now. Then by (6.6) and Lemma 7, there is $\{\tau_n\} \subset \{t'_n\}$, and $\boldsymbol{\lambda} \in \mathbb{C}^V$, such that $\lim_{n\to\infty} \tau_n^{-1}\boldsymbol{\Delta}_{\tau_n} = \boldsymbol{\lambda}$. Divide both sides of (6.4) and (6.3) by $\tau_n$ and let $n \to \infty$. Then there is $z \in \mathbb{C}$, such that $\boldsymbol{\lambda} = z\mathbf{u}$.

To find the value of $z$, divide both ends of (6.5) by $t^2$ and let $t \to 0$ through $\tau_n$. By Lemma 2,

$$iv_j(\mathbf{c}\cdot\mathbf{q})u_j z = \tfrac{1}{2}H(z\mathbf{u}) - v_j\left(-iK + \tfrac{1}{2}(\mathbf{c}\cdot\mathbf{q})^2 - \tfrac{1}{2}\sum_{\mathbf{n}\in\mathcal{V}^*} p_j(\mathbf{n})c_\mathbf{n}^2\right)$$

$$+ iv_j \sum_{\mathbf{n}\in\mathcal{V}^*} p_j(\mathbf{n})c_\mathbf{n}(\mathbf{n}\cdot\mathbf{u})z.$$

It is then routine to check that $z$ is the solution to (6.1) with negative real part. The remaining part of the proof is similar to that for Theorem 2 and thus is omitted.

PROOF OF LEMMA 7. If $\liminf_{n\to\infty} |t'_n/\Delta_{k_0}(t'_n)| = 0$, then choose $\{\tau_n\} \subset \{t'_n\}$ such that
$$\lim_{n\to\infty} \frac{\tau_n}{\Delta_{k_0}(\tau_n)} = 0.$$

Divide both sides of (6.3) and (6.4) by $\Delta_{k_0}(t)$ and let $t \to 0$ through $\tau_n$. Then (6.6) leads to $\boldsymbol{\xi} = M\boldsymbol{\xi}$. Because $\xi_{k_0} = 1$, $\boldsymbol{\xi} = u_{k_0}^{-1}\mathbf{u}$. Now divide (6.5) by $\Delta_{k_0}(t)^2$ and let $t \to 0$ through $\tau_n$. Then it is seen that $H(\boldsymbol{\xi}) = 0$, implying $H(\mathbf{u}) = 0$, which is a contradiction. $\square$

**7. An estimator for the right eigenvector of the mean matrix.** This section gives the proof for Theorem 4. We need a few lemmas.

LEMMA 8. *Assume the same conditions as in Theorem 4. For $\lambda > 0$, denote*
$$S_j(\lambda) = E_j S(\omega, \lambda), \qquad j \in \mathcal{V}, \qquad \mathbf{S}_\lambda = (S_1(\lambda), \ldots, S_V(\lambda)).$$
*Then for $\lambda \in (0,1)$, $S_j(\lambda) < \infty$ and*

(7.1) $$\mathbf{S}_\lambda^t = (1 - \lambda M)^{-1}\mathbf{1}^t,$$

(7.2) $$E_j(S(\omega,\lambda))^t = O\left(\frac{1}{(1-\lambda)^{2t-1}}\right) \quad \text{as } \lambda \to 1, \ t = 2, 3, 4, \ j \in \mathcal{V}.$$

LEMMA 9.

(7.3) $$\lim_{\lambda \nearrow 1}(1-\lambda)\mathbf{S}_\lambda = \mathbf{u}.$$

Assume the lemmas to be true for now. By Lemma 8, $S(\omega, \lambda)$ is integrable. Letting
$$\widetilde{S}(\omega, \lambda) = S(\omega, \lambda) - S_k(\lambda),$$
by Lemma 9, we need to show that if $\omega_1, \omega_2, \ldots$ are i.i.d. $\sim P_k$, then
$$\frac{1-\lambda_N}{N} \sum_{n=1}^{N} \widetilde{S}(\omega_n, \lambda_N) \to 0, \qquad P_k\text{-a.s.}$$

To this end, by Borel–Cantelli and the Markov inequality, it is enough to show that
$$\sum_{N=1}^{\infty} \frac{(1-\lambda_N)^4}{N^4} E_k\left(\sum_{n=1}^{N} \widetilde{S}(\omega_n, \lambda_N)\right)^4 < \infty.$$

Because $\omega_n$ are i.i.d., and $E_k[\widetilde{S}(\omega, \lambda_N)] = 0$, by (7.2),

$$E_k\left(\sum_{n=1}^{N} \widetilde{S}(\omega_n, \lambda_N)\right)^4 = 3N(N-1)[\mathrm{Var}_k S(\omega, \lambda_N)]^2 + N E_k \overline{S}_k(\omega, \lambda_N)^4$$

$$\leq \frac{CN^2}{(1-\lambda_N)^6} + \frac{CN}{(1-\lambda_N)^7}, \qquad N \to \infty$$

for some constant $C$. Therefore, by (2.16),

$$\sum_{N=1}^{\infty} \frac{(1-\lambda_N)^4}{N^4} E_k \left(\sum_{n=1}^{N} \widetilde{S}(\omega_n, \lambda_N)\right)^4$$

$$\leq \sum_{N=1}^{\infty} \frac{C}{(1-\lambda_N)^2 N^2}\left(1 + \frac{1}{(1-\lambda_N)N}\right) < \infty.$$

PROOF OF LEMMA 8. That $S_j(\lambda) < \infty$ for all $j \in \mathcal{V}$ and $\lambda \in (0,1)$ and (7.1) are easy consequences of recursion. For $\lambda \in (0,1)$, because $\mathbf{v}(1 - \lambda M)^{-1} = (1-\lambda)^{-1}\mathbf{v}$,

$$(7.4) \qquad \mathbf{v} \cdot \mathbf{S}_\lambda = \mathbf{v}(1 - \lambda M)^{-1}\mathbf{1}^t = \frac{\mathbf{v} \cdot \mathbf{1}}{1-\lambda} = \frac{1}{1-\lambda}.$$

Because all components of $\mathbf{v}$ are strictly positive, $S_j(\lambda) \sim (1-\lambda)^{-1}$, $\lambda \to 1-$, and thus (7.2) is proved for $t = 1$. The proof of (7.2) is similar for $t = 2$, 3 and 4. We shall show the details of the proof for $t = 4$, to illustrate how the indices of $1 - \lambda$ in the asymptotics (7.2) are counted.

Suppose we have shown $E_j S(\omega, \lambda)^t < \infty$ and (7.2) for $t = 2, 3$. Given a tree $\omega$ rooted with $j$, suppose the branching rule applied by the root is $(j \to \mathbf{n})$. For $s \in \mathcal{V}$, and $l = 1, \ldots, n_s$, let $\omega_{l,s}$ be the subtree rooted with the $l$th particle of type $s$ in $\mathbf{n}$. For each $D \geq 0$, define $S_D(\omega, \lambda) = \sum_{x \in \omega} \lambda^{|x|}\mathbf{1}_{\{|x| \leq D\}}$. Then

$$S_D(\omega, \lambda) = 1 + \Sigma(\lambda) \stackrel{\Delta}{=} 1 + \lambda \sum_{s \in \mathcal{V}} \sum_{l=1}^{n_s} S_{D-1}(\omega_{l,s}, \lambda).$$

Then

$$(7.5) \quad E_j S_D(\omega, \lambda)^4 = 1 + 4E_j\Sigma(\lambda) + 6E_j\Sigma^2(\lambda) + 4E_j\Sigma^3(\lambda) + E_j\Sigma^4(\lambda).$$

For $t = 1, 2, 3, 4$, and $\mathbf{n} \in \mathcal{V}^*$, let

$m_t(\lambda) = \max\{E_s S(\omega, \lambda)^t : s \in \mathcal{V}\}$,

$I_t(\mathbf{n}) = \{((l_1, s_1), \ldots, (l_t, s_t)) : 1 \leq l_i \leq n_{s_i}, \ i = 1, \ldots, t,$

$\qquad\qquad$ and $(l_i, s_i)$ are different from each other$\}$.

Then by the multinomial expansion, it is seen that

$$E_j[\Sigma^4(\lambda)] = \lambda^4 \sum_{s \in \mathcal{V}} M_{js} E_s[S_{D-1}(\omega, \lambda)^4] + \lambda^4 \sum_{\mathbf{n} \in \mathcal{V}^*} p_j(\mathbf{n}) \sum_{i=1}^{4} \Sigma_{n,i}(\lambda),$$

where, as $\lambda \to 1-$,

$$\Sigma_{n,1}(\lambda) = \sum_{I_2(\mathbf{n})} E_j[S_{D-1}(\omega_{l_1,s_1}, \lambda)^3] E_j[S_{D-1}(\omega_{l_2,s_2}, \lambda)]$$

$$\leq |I_2(\mathbf{n})| m_1(\lambda) m_3(\lambda) = O(|\mathbf{n}|^4 (1-\lambda)^{-6}),$$

with the summation over all $((l_1, s_1), (l_2, s_2)) \in I_2(\mathbf{n})$, and likewise,

$$\Sigma_{n,2}(\lambda) = \sum_{I_2(\mathbf{n})} E_{l_1}[S_{D-1}(\omega_{l_1,s_1}, \lambda)^2] E_{l_2}[S_{D-1}(\omega_{l_2,s_2}, \lambda)^2]$$

$$\leq |I_2(\mathbf{n})| m_2^2(\lambda) = O(|\mathbf{n}|^4 (1-\lambda)^{-6}),$$

$$\Sigma_{n,3}(\lambda) = \sum_{I_3(\mathbf{n})} E_j[S_{D-1}(\omega_{l_1,s_1}, \lambda)^2] \prod_{i=2}^{3} E_j[S_{D-1}(\omega_{l_i,s_i}, \lambda)]$$

$$\leq |I_3(\mathbf{n})| m_2(\lambda) m_1(\lambda)^2 = O(|\mathbf{n}|^4 (1-\lambda)^{-5}),$$

$$\Sigma_{n,4}(\lambda) = \sum_{I_4(\mathbf{n})} \prod_{i=1}^{4} E_j[S_{D-1}(\omega_{l_i,s_i}, \lambda)]$$

$$\leq |I_4(\mathbf{n})| m_1(\lambda)^4 = O(|\mathbf{n}|^4 (1-\lambda)^{-4}).$$

By (2.15), $\sum_{\mathbf{n} \in \mathcal{V}^*} p_j(\mathbf{n}) |\mathbf{n}|^4 < \infty$. Therefore, the above estimates imply that there is a constant $C$, such that as $\lambda \to 1-$,

$$E_j[\Sigma^4(\lambda)] \leq \lambda^4 \sum_{l \in \mathcal{V}} M_{jl} E_l[S_{D-1}(\omega, \lambda)^4] + C(1-\lambda)^{-6}.$$

By Hölder's inequality, the other summands on the right-hand side of (7.5) are dominated by $(1-\lambda)^{-6}$. Therefore, for some constant, still denoted $C$,

$$E_j[S_D(\omega, \lambda)^4] \leq \lambda^4 \sum_{s \in \mathcal{V}} M_{js} E_s[S_{D-1}(\omega, \lambda)^4] + C(1-\lambda)^{-6}.$$

Multiply both sides of the above inequality by $v_j$ to the left, and sum over $j \in \mathcal{V}$. Let

$$A_D = \sum_{k \in \mathcal{V}} v_k E_k[S_D(\omega, \lambda)^4].$$

Then

$$A_D \leq \lambda^4 A_{D-1} + C(1-\lambda)^{-6} \implies \lim_{D \to \infty} A_D \leq \frac{C}{(1-\lambda)^{-7}}.$$

Since $v_k \in (0,1)$, $k \in \mathcal{V}$, the last formula leads to (7.2). □

PROOF OF LEMMA 9. Let $\boldsymbol{\xi}_\lambda = (1-\lambda)\mathbf{S}_\lambda$. Then by (7.4),

$$\mathbf{v} \cdot \boldsymbol{\xi}_\lambda = 1.$$

Note that all the coordinates of $\boldsymbol{\xi}_\lambda = (1-\lambda)\mathbf{S}_\lambda$ are positive. Indeed, by (7.1),

$$\mathbf{S}_\lambda^t = \sum_{n=0}^\infty \lambda^n M^n \mathbf{1}^t.$$

Every $M^n$ is nonnegative. Thus all the coordinates of $\mathbf{S}_\lambda$, hence all those of $\boldsymbol{\xi}_\lambda$, are positive. Since all the coordinates of $\mathbf{v}$ are strictly positive, then $\boldsymbol{\xi}_\lambda$ is bounded and thus has cluster points. If $\boldsymbol{\xi}$ is a cluster point of $\boldsymbol{\xi}_\lambda$ as $\lambda \to 1-$, then $\mathbf{v} \cdot \boldsymbol{\xi} = 1$. On the other hand,

$$(1-M)\boldsymbol{\xi} = \lim_{\lambda \to 1-} (1-\lambda M)\boldsymbol{\xi}$$
$$= \lim_{\lambda \to 1-} (1-\lambda M)(1-\lambda)(1-\lambda M)^{-1}\mathbf{1} = \lim_{\lambda \to 1-} (1-\lambda)\mathbf{1} = 0$$

and hence $\boldsymbol{\xi} = M\boldsymbol{\xi}$. Thus $\boldsymbol{\xi} = \mathbf{u}$. □

**Acknowledgment.** The author thanks an anonymous referee for kindly pointing out several important relevant articles.

DEPARTMENT OF STATISTICS
UNIVERSITY OF CHICAGO
CHICAGO, ILLINOIS 60637
USA
E-MAIL: chi@galton.uchicago.edu